\documentclass[11pt,authoryear,preprint]{elsarticle}
\usepackage{amsmath}
\usepackage{amssymb}
\usepackage{bm}
\usepackage{mathrsfs,dsfont}
\usepackage{graphicx,color}
\usepackage[footnotesize]{caption}
\usepackage[scriptsize]{subfigure}
\RequirePackage[colorlinks,citecolor=blue,urlcolor=blue]{hyperref}
\RequirePackage{natbib}
\usepackage{theorem}
\usepackage{xspace}
\usepackage{algorithm}
\usepackage{algorithmic}
\usepackage[english]{babel}

\allowdisplaybreaks

\theoremstyle{definition}

\newtheorem{proposition}{\sc{Proposition}}
\newtheorem{corollary}{\sc{Corollary}}

\newcommand{\indic}{\mathds{1}}

\newcommand{\Y}{\mathbb{Y}}

\newcommand\Be{\text{Be}\xspace}
\newcommand\truncNorm{\text{t}\mathcal{N}_{[0,1]}\xspace}

\newcommand{\R}{\mathds{R}}
\newcommand{\E}{\mathds{E}}
\newcommand{\N}{\mathds{N}}
\newcommand{\Nbb}{\mathds{N}}
\renewcommand{\P}{\mathds{P}}

\newcommand{\ddr}{\mathrm{d}}
\newcommand{\edr}{\mathrm{e}}

\newcommand{\Lc}{\mathcal{L}}

\newcommand{\ep}{\varepsilon}

\def\simiid{\stackrel{\mbox{\scriptsize{iid}}}{\sim}}

\newcommand{\mt}{\tilde{\mu}} % misura completamente aleatoria (CRM)

\newcommand{\St}{\tilde{S}}
\makeatletter
\def\ps@pprintTitle{%
 \let\@oddhead\@empty
 \let\@evenhead\@empty
 \def\@oddfoot{}%
 \let\@evenfoot\@oddfoot}
\makeatother
\begin{document}

\begin{frontmatter}

\title{\textbf{Full Bayesian inference\\with hazard mixture models}}

\author[cca]{Julyan Arbel}\ead{julyan.arbel@carloalberto.org}
\author[pavia,cca]{Antonio Lijoi}\ead{lijoi@unipv.it}
\author[torino,cca]{Bernardo Nipoti\corref{cor1}}\ead{bernardo.nipoti@carloalberto.org}

\address[cca]{Collegio Carlo Alberto, via Real Collegio, 30, 10024 Moncalieri, Italy}
\address[pavia]{Department of Economics and Management, University of Pavia, Via San Felice 5, 27100 Pavia, Italy}
\address[torino]{Department of Economics and Statistics, University of Torino, C.so Unione Sovietica 218/bis, 10134 Torino, Italy}
\cortext[cor1]{Corresponding author. Phone: +39-011-6705023}

\begin{abstract}
Bayesian nonparametric inferential procedures based on Markov chain Monte Carlo marginal methods typically yield point estimates in the form of posterior expectations. Though very useful and easy to implement in a variety of statistical 
problems, these methods may suffer from some limitations if used to estimate non-linear functionals of the posterior distribution. The main goal of the present paper is to develop a novel methodology that extends a well-established 
marginal procedure designed for hazard mixture models, in order to draw approximate inference on survival functions that is not limited to the posterior mean but includes, as remarkable examples, credible intervals and median survival time.
Our approach relies on a characterization of the posterior moments that, in turn, is used to approximate the posterior distribution by means of a technique based on Jacobi polynomials. The inferential performance of our methodology is analyzed by means of an extensive study of simulated data and real data consisting of leukemia remission times.  Although tailored to the survival analysis context, the procedure we introduce can be adapted to a range of other models for which moments of the posterior distribution can be estimated.
\end{abstract}


\end{frontmatter}

\section{Introduction}
Most commonly used inferential procedures in Bayesian nonparametric practice rely on the implementation of sampling algorithms that can be gathered under the general umbrella of Blackwell--MacQueen P\'olya urn schemes. These are characterized by the marginalization with respect to an infinite-dimensional random element that defines the de Finetti measure of an exchangeable sequence of observations or latent variables. Henceforth we will refer to them as \emph{marginal methods.} Besides being useful for the identification of the basic building blocks of ready to use Markov chain Monte Carlo (MCMC) sampling strategies, marginal methods have proved to be effective for an approximate evaluation of Bayesian point estimators in the form of posterior means. They are typically used with models for which the predictive distribution is available in closed form. Popular examples are offered by mixtures of the Dirichlet process for density estimation \citep{EscWes95} and mixtures of gamma processes for hazard rate estimation \citep{IshJam04}. While becoming well-established tools, these computational techniques are easily accessible also to practitioners through a straightforward software implementation \citep[see for instance][]{jara2011dppackage}. Though it is important to stress their relevance both in theory and in practice, it is also worth pointing out that 
Blackwell--MacQueen P\'olya urn schemes suffer from some drawbacks which we wish to address here. Indeed, one easily notes that the posterior estimates provided by marginal methods are not suitably endowed with measures of uncertainty such as posterior credible intervals. Furthermore, using the posterior mean as an estimator is equivalent to choosing a square loss function whereas in many situations of interest other choices such as absolute error or $0\mbox{--}1$ loss functions and, as corresponding estimators, median or mode of the posterior distribution 
of the survival function, at any fixed time point $t$, would be preferable. Finally, they do not naturally allow inference on functionals of the distribution of survival times, such as the median survival time, to be drawn. 
A nice discussion of  these issues is provided by \cite{GelKot02} where the focus is on mixtures of the Dirichlet process: the authors suggest complementing the use of marginal methods with a sampling strategy that aims at generating approximate trajectories of the Dirichlet process from its truncated stick-breaking representation. 

The present paper aims at proposing a new procedure that combines closed-form analytical results arising from the application of marginal methods with an approximation of the posterior distribution which makes use of posterior moments. The whole machinery is developed for the estimation of survival functions that are modeled in terms of hazard rate functions. To this end, let $F$ denote the cumulative distribution function (CDF) associated to a probability distribution on $\R^+$. The corresponding survival and cumulative hazard functions are denoted as \
\[
S(t)=1-F(t)
\quad\mbox{ and } \quad
H(t)=-\int_{[0,t]}\frac{\ddr F(s)}{F(s-)},
\] 
for any $t>0$, respectively, where $F(s-):=\lim_{\ep\downarrow 0}F(s-\ep)$ for any positive $s$. 
If $F$ is absolutely continuous, one has $H(t)=-\log (S(t))$ and the hazard rate function associated to $F$ is, thus, defined as $h(t)=F'(t)/[1-F(t-)]$. It should be recalled that survival analysis has been one of the most relevant areas of application of Bayesian nonparametric methodology soon after the groundbreaking contribution of  \cite{Ferguson73}. A number of papers in the '70s and the '80s have been devoted to the proposal of new classes of priors that accommodate for a rigorous analytical treatment of Bayesian inferential problems with censored survival data. Among these we need to mention neutral to the right processes proposed in \cite{Doksum} and used to define a prior for the CDF $F$: since they share a conjugacy property they represent a tractable tool for drawing posterior inferences. 
Another noteworthy class of priors has been proposed in \cite{Cervo}, where a beta process is used as a  nonparametric prior for the cumulative hazard function $H$ has been proposed. Also in this case, one can considerably benefit from a useful conjugacy property.
  
As already mentioned, we plan to propose a method for full Bayesian analysis of survival data by specifying a prior on the hazard rate $h$. The most popular example is the gamma process mixture that has been originally proposed in \cite{DykLau81} and generalized in later work by \cite{LoWen89} and \cite{Jam05} to include any mixing random measure and any mixed kernel. Recently \cite{LijNip14} have extended such framework to the context of partially exchangeable observations. The uses of random hazard mixtures in practical applications 
have been boosted by the recent developments of powerful computational techniques that allow for an approximate evaluation of posterior inferences on quantities of statistical interest. Most of these arise from a marginalization with respect to a completely random measure that identifies the de~Finetti measure of the exchangeable sequence of observations. 
See, e.g., \cite{IshJam04}. Though they are quite simple to implement, the direct use of their output can only yield point estimation of the hazard rates, or of the survival functions, at fixed time points through posterior means. The main goal of the present paper is to show that a clever use of a moment-based approximation method does provide a relevant upgrade on the type of inference one can draw via marginal sampling schemes. The takeaway message is that the information gathered by marginal methods is not confined to the posterior mean but is actually much richer and, if properly exploited, can lead to a more complete posterior inference. To understand this, we shall refer to a sequence of exchangeable survival times $(X_i)_{i\ge 1}$ such that $\P[X_1>t_1,\ldots,X_n>t_n\,|\,\tilde P]=\prod_{i=1}^n \tilde{S}(t_i)$ where $\tilde P$ is a random probability measure on $\R^+$ and $\tilde{S}(t)=\tilde P((t,\infty))$ is the corresponding random survival function. Given a suitable sequence of latent variables $(Y_i)_{i\ge 1}$, we will provide a closed-form expression for 
\begin{equation}
\label{eq:moments_r}
\E[\tilde{S}^r(t)\,|\,\bm{X},\bm{Y}], \qquad\mbox{ for any } r\ge 1,\:\:\mbox{ and }t>0,
\end{equation}
with $\bm{X}=(X_1,\ldots,X_n)$ and $\bm{Y}=(Y_1,\ldots,Y_n)$. Our strategy consists in approximating the posterior distribution of $\tilde S(t)$, at each instant $t$, and relies on the fact that, along with the posterior mean, marginal models allow to straightforwardly estimate posterior moments of any order of $\tilde S(t)$. Indeed, an MCMC sampler yields a sample from the posterior distribution of $\bm{Y}$ given $\bm{X}$: this can be used to integrate out the latent variables appearing in \eqref{eq:moments_r} and obtain a numerical approximate evaluation of the posterior moments $\E[\tilde{S}^r(t)\,|\,\bm{X}]$. 
These are finally used to deduce, with almost negligible effort, an approximation of the posterior distribution of $\tilde S(t)$ and, in turn, to estimate its functionals.

It is to be mentioned that one could alternatively resort to a different approach that boils down to the simulation of the trajectories of the completely random measure that defines the underlying random probability measure from its posterior distribution. In density estimation problems, this is effectively illustrated in \cite{NPW04}, \cite{NP09} and \cite{BLNP13}. As for hazard rates mixtures estimation problems, one can refer to \cite{Jam05}, \cite{NieWal04} and \cite{Nie13}. In particular, \cite{Jam05} provides a posterior characterization that is the key for devising a \cite{FerKla72} representation of the posterior distribution of the completely random measure which enters the definition of the prior for the hazards. Some numerical aspects related to the implementation of the algorithm can be quite tricky since one needs to invert the L\'evy intensity to simulate posterior jumps and a set of suitable latent variables need to be introduced in order to sample from the full conditionals of the hyperparameters. These aspects are well described and addressed in \cite{Nie13}.
 
The paper is organized as follows. In Section~\ref{sec:haz} we briefly review hazard mixture models and recall some of their most important properties. 
We further provide explicit expressions characterizing the posterior moments of any order of a random survival function, both for general framework and for the extended gamma process case. Section~\ref{sec:approx} is dedicated to the problem of approximating the distribution of a random variable on $[0,1]$, provided that the first $N$ moments are known. In particular, in Section~\ref{sec:approx1} we describe a convenient methodology based on Jacobi polynomials, whereas in Section~\ref{sec:approx2} we apply such methodology in order to approximate random survival functions and we perform a numerical study to test its performance. In Section~\ref{sec:real} we focus on the use of the introduced methodology for carrying out Bayesian inference on survival functions. Specifically, in Section~\ref{algo} we present the algorithm, whereas in Sections~\ref{sec:simulated} and \ref{sec:leukemia} we analyze, respectively, simulated data a real two-sample dataset on leukemia remission times. For the sake of exposition simplicity, we postponed to the Appendix technicalities such as expressions for the full conditional distributions involved in the algorithm and instructions on how to take into account the presence of censored data.

\section{Hazard mixture models}\label{sec:haz}
A well-known nonparametric prior for the hazard rate function within multiplicative intensity models used in survival analysis arises as a mixture of \textit{completely random measures} (CRMs). To this end, recall that a CRM $\mt$ on a space $\Y$ is a boundedly finite random measure that, when evaluated at any collection of pairwise disjoint sets $A_1,\ldots,A_d$, gives rise to mutually independent random variables $\mt(A_1),\ldots,\mt(A_d)$, for any $d\ge 1$. Importantly, CRMs are almost surely discrete measures \citep{Kin92}. A detailed treatment on CRMs can also be found in \cite{DalVer03}. With reference to Theorem~1 in \cite{Kingman67}, we shall assume that $\mt$ has no fixed atoms, which in turn implies the existence of a measure $\nu$ on $\R^+\times \Y$ such that
$ \int_{\R^+\times\Y}\min\{s,1\}\nu(\ddr s, \ddr y)<\infty$
and
\begin{equation}
\label{eq:lapl_univ}
 \E\left[\edr^{-\int_{\Y}f(y)\mt(\ddr y)}\right]=\exp\left(-\int_{\R^+\times\Y}
\left[1-\exp\left(-s\,f(y)\right)\right]\nu(\ddr s, \ddr y)\right),
\end{equation}
for any measurable function $f:\Y\rightarrow \R$ such that $\int_\Y\left|f\right|\,\ddr\mt<\infty$, with probability~1. The measure $\nu$ is termed the \textit{L\'evy intensity} of $\mt$. For our purposes, it will be useful to rewrite $\nu$ as
\begin{equation*}
\label{eq:levy_intens}
 \nu(\ddr s,\ddr y)=\rho_y(s)\,\ddr s\,c\,P_0(\ddr y),
\end{equation*}
where $P_0$ is a probability measure on $\Y$, $c$ a positive parameter, and $\rho_y(s)$ is some transition kernel on  $\Y\times\R^+$. If $\rho_y=\rho$, for any $y$ in $\Y$, the CRM $\mt$ is said \textit{homogeneous}. Henceforth, we further suppose that $P_0$ is non-atomic. A well-known example corresponds to $\rho_y(s)=\rho(s)=\edr^{-s}/s$, for any $y$ in $\Y$, which identifies a so-called \textit{gamma CRM}. With such a choice of the L\'evy intensity,  it can be seen, from \eqref{eq:lapl_univ}, that for any $A$ %in $\Yf$ 
such that $P_0(A)>0$, the random variable $\mt(A)$ is gamma distributed, with shape parameter 1 and rate parameter $cP_0(A)$.
If $k(\,\cdot\,;\cdot\,)$ is a transition kernel on $\R^+\times \Y$, a prior for $h$ is the distribution of the random hazard rate (RHR)  
\begin{equation}\label{RHR}
\tilde h(t)=\int_\Y k(t;y)\mt(\ddr y),
\end{equation}
where $\mt$ is a CRM on $\Y$. We observe that, if $
\lim_{t\rightarrow \infty} \int_0^t \tilde h(s)\ddr s=\infty %\qquad \P\mbox{-a.s.},
$
with probability 1, then one can adopt the following model 
\begin{equation}
\label{eq:exch_model}
  \begin{split}
    X_i\,|\,\tilde P \: &\simiid \: \tilde P\\
    \tilde P((\,\cdot\,,\infty)) \:
    &\stackrel{\scriptsize{\mbox{d}}}{=}\: \exp\left(-\int_0^{\,\cdot}
      \tilde h(s)\,\ddr s\right)
  \end{split}
\end{equation}
for a sequence of (possibly censored) survival data $(X_i)_{i\ge 1}$. This means that $\tilde h$ in (\ref{RHR}) defines a random  survival function 
$
t\mapsto\tilde S(t)=\exp(-\int_0^t \tilde h(s)\ddr s).$
In this setting, \cite{DykLau81} characterize the posterior distribution of the so-called \textit{extended gamma process}: this is obtained when $\mt$ is a gamma CRM and $k(t;y)=\indic_{(0,t]}(y)\,\beta(y)$ for some positive right-continuous function $\beta :\R^+\to\R^+$. The same kind of result is proved in \cite{LoWen89} for \textit{weighted gamma processes} corresponding to RHRs obtained when $\mt$ is still a  gamma CRM and $k(\,\cdot\,;\,\cdot\,)$ is an arbitrary kernel. Finally, a posterior characterization has been derived in \cite{Jam05} for any CRM $\mt$ and kernel $k(\,\cdot\,;\,\cdot\,)$. 
\bigskip\\
We shall quickly display such a characterization since it represents the basic result our construction relies on. For the ease of exposition we confine ourselves to the case where all the observations are exact, the extension to the case that includes right-censored data being straightforward and detailed in \cite{Jam05}. See also \ref{sec:censored}. For an $n$-sample $\bm{X}=(X_1,\ldots,X_n)$ of exact data, the likelihood function equals
\begin{equation}\label{eq:likelihood}
 \Lc(\mt;\bm{X})=e^{-\int_\Y K_{\bm{X}}(y)\mt(\ddr y)}\prod_{i=1}^n\int_\Y k(X_i;y)\mt(\ddr y),
\end{equation}
where $K_t(y) = \int_0^{t}k(s;y)\ddr s$ and $K_{\bm{X}}(y)=\sum_{i=1}^n K_{X_i}(y)$. A useful augmentation suggests introducing latent random variables $\bm{Y}=(Y_1,\ldots,Y_n)$ such that the joint distribution of $(\mt,\bm{X},\bm{Y})$ coincides with
\begin{equation}\label{joo}
 %\Lc(\mt,\bm{X},\bm{Y})=
 e^{-\int_\Y K_{\bm{X}}(y)\mt(\ddr y)}\prod_{i=1}^n k(X_i;Y_i)\mt(\ddr Y_i)\, Q(\ddr \mt),
\end{equation}
where $Q$ is the probability distribution of the completely random measure $\mt$, characterized by the Laplace transform functional in \eqref{eq:lapl_univ} \citep[see for instance][]{DalVer03}. 
The almost sure discreteness of $\tilde \mu$ implies there might be ties among the $Y_i$'s with positive probability.
Therefore, we denote the distinct values among $\bm{Y}$ with $(Y_1^*,\ldots,Y_k^*)$, where $k\leq n$,  and, for any $j=1,\ldots,k$, we define $C_j=\left\{l\,:\,Y_l=Y_j^*\right\}$ and $n_j=\#C_j$, the cardinality of $C_j$. Thus, we may rewrite the joint distribution in \eqref{joo} as
\begin{equation}
\label{eq:joint.rhr.1}
%\Lc(\mt,\bm{X},\bm{Y})=
e^{-\int_\Y K_{\bm{X}}(y)\mt(\ddr y)}\prod_{j=1}^k \mt(\ddr Y_j^*)^{n_j}
\prod_{i\in C_j}k(X_i;Y_j^*)\, Q(\ddr \mt). 
\end{equation}
We introduce, also, the density function
\begin{equation}
  \label{eq:jump_dens}
  f(s\,|\,\kappa,\xi,y)\propto s^\kappa\,\edr^{-\xi s}\,\rho_y(s)\:\indic_{\R^+}(s)
\end{equation}
for any $\kappa\in \N$ and $\xi>0$. The representation displayed in \eqref{eq:joint.rhr.1}, combined with results concerning disintegrations of Poisson random measures, leads to prove the following

\begin{proposition}%[\citealp{Jam05}]
\label{Jam}{\rm (\citealp{Jam05})}
Let $\tilde h$ be a RHR as defined in \eqref{RHR}. The posterior distribution of $\tilde h$, given $\bm{X}$ and $\bm{Y}$, coincides with the distribution of the random hazard
\begin{equation}\label{post}
 %\mt^*\stackrel{\mbox{\scriptsize{d}}}{=}
  \tilde{h}^*+\sum_{j=1}^k J_j k(\,\cdot\,;Y_j^*),
\end{equation}
where $\tilde{h}^*(\,\cdot\,)=\int_\Y k(\,\cdot\,;y)\,\mt^*(\ddr y)$ and $\mt^*$ is a CRM without fixed points of discontinuity whose L\'evy intensity is
\begin{equation*}
 \nu^*(\ddr s, \ddr y)= e^{-s K_{\bm{X}}(y)}\rho_y (s)\ddr s\,c P_0(\ddr y).
\end{equation*}
The jumps $J_1,\ldots,J_k$ are mutually independent and independent of $\mt^*$. Moreover, for every $j=1,\ldots,k$, the distribution of the jump $J_j$ has density function $f(\,\cdot\,|\, n_j,K_{\bm{X}}(Y_j^*),Y_j^*)$ with $f$ defined in \eqref{eq:jump_dens}.
\end{proposition}
\medskip

See \cite{LPW08} for an alternative proof of this result. The posterior distribution of $\tilde h$ displays a structure that is common to models based on CRMs, since it consists of the combination of two components: one without fixed discontinuities and the other with jumps at fixed points. In this case, the points at which jumps occur coincide with the distinct values of the latent variables $Y_1^*,\ldots,Y_k^*$. Furthermore, the distribution of the jumps $J_j$ depends on the respective locations $Y_j^*$.

Beside allowing us to gain insight on the posterior distribution of $\tilde h$, Proposition~\ref{Jam} is also very convenient for simulation purposes. See, e.g., \cite{IshJam04}. Indeed, \eqref{post} allows obtaining an explicit expression for the posterior expected value of $\tilde S(t)$ (or, equivalently, of $\tilde h(t)$), for any $t>0$,  conditionally on the latent variables $\bm{Y}$. One can, thus, integrate out the vector of latent variables $\bm{Y}$, by means of a Gibbs type algorithm, in order to approximately evaluate the posterior mean of $\tilde S(t)$ (or $\tilde h(t)$). As pointed out in next section, a combination of Proposition~\ref{Jam} and of the same Gibbs sampler we have briefly introduced actually allows moments of $\tilde{S}(t)$, of any order, to be estimated. We will make use of the first $N$ of these estimated moments to approximate, for each $t>0$, the posterior distribution of $\tilde{S}(t)$ and therefore to have the tools for drawing meaningful Bayesian inference. The choice of a suitable value for  $N$ will be discussed in Section~\ref{sec:approx2}.

As pointed out in the Introduction, one can, in line of principle, combine Proposition~\ref{Jam} with the Ferguson and Klass representation to undertake an alternative approach that aims at simulating the trajectories from the posterior distribution of the survival function. This can be achieved by means of a Gibbs type algorithm that involves sampling $\tilde \mu^*$ and $Y^*_j$, for $j=1,\ldots,k$, from the corresponding full conditional distributions. Starting from the simulated trajectories one could then approximately evaluate all the posterior quantities of interest. Since this approach does not rely on the marginalization with respect to the mixing CRM $\tilde\mu$, we refer to it as an example of non-marginal, or \textit{conditional}, method. 
An illustration, with an application to survival analysis, is provided in \cite{Nie13} and it appears that the approach, though achievable, may be difficult to implement. The main non-trivial issues one has to deal with are the inversion of the L\'evy measure, needed to sample the jumps, and the sampling from the full conditionals of the hyperparameters. The latter has been addressed by  \cite{Nie13} through a clever augmentation scheme that relies on a suitable collection of latent variables. In any case, it is worth recalling that even the Ferguson and Klass algorithm is based on an approximation since a realization of $\tilde \mu^*$ is approximated with a finite number of jumps.

In the next sections we will focus on marginal methods since our goal is to show that they allow for a full Bayesian inference, beyond the usual evaluation of posterior means. The required additional effort to accomplish this task is minimal and boils down to computing a finite number of posterior moments of $\tilde S(t)$, at a given $t$. An approximate evaluation of these moments can be determined by resorting to 
\eqref{post} which yields closed-form expressions for the posterior moments of the random variable $\St(t)$, conditionally on both the data $\bm{X}$ and the latent variables $\bm{Y}$. 
\begin{proposition}\label{mom} For every $t>0$ and $r>0$,
\begin{multline*}
\E[\St^r(t)\,|\,\bm{X},\bm{Y}]=\exp\left\{-c\int_{\R^+\times \Y}\left(1-\edr^{-r K_t(y) s}\right)\edr^{-K_{\bm{X}}(y) s} \rho(s) \ddr s P_0(\ddr y)\right\}\\
\times\prod_{j=1}^k \frac{1}{\mbox{B}_j}\int_{\R^+}\exp\left\{-s\left(r K_t(Y_j^*)+K_{\bm{X}}(Y_j^*)\right)\right\}s^{n_j}\rho(s)\ddr s,
\end{multline*}
where $B_j=\int_{\R^+} s^{n_j} \exp\left\{-s K_{\bm{X}}(Y_j^*)\right\}\rho(s) \ddr s$, for $j=1,\ldots,k$.\\
\end{proposition} 
Although the techniques we will describe in next section hold true for any specification of $\mt$ and kernel $k(\,\cdot\,;\cdot\,)$, for illustration purposes we focus on the extended gamma process case \citep{DykLau81}. 
More specifically, we consider a kernel  $k(t;y)=\indic_{(0,t]}(y)\beta$, with $\beta>0$. 
This choice of kernel is known to be suitable for modeling monotone increasing hazard rates and to give rise to a class of random hazard functions with nice asymptotic properties \citep{DeBPecPru09}. 
Moreover, without loss of generality, we suppose that $X_1>X_2>\ldots>X_n$, we set, for notational convenience, $X_0\equiv\infty$ and $X_{n+1}\equiv 0$ and we introduce $\xi_l\equiv \sum_{i=1}^l X_i$, for any $l\geq 1$, and set $\xi_0\equiv 0$. The next Corollary displays an expression for the conditional moments corresponding to this prior specification.
\begin{corollary} \label{cor:moments}
For every $t>0$ and $r>0$,
\begin{multline}\label{mom_ext_gamma}
\E[\St^r(t)\,|\,\bm{X},\bm{Y}]=
%\E[\St(t)^r\,|\,\mathbf{X},\mathbf{Y}]=
\prod_{i=0}^{n}\exp\left\{-c\int_{X_{i+1}\wedge t}^{X_{i}\wedge t}\log\left(1+r\,\frac{t-y}{\xi_{i}-iy+1/\beta}\right) P_0 (\ddr y)\right\}\\
\times\prod_{j=1}^k \left(1+r\,\frac{(t-Y_j^*)\indic_{[Y_j^*,\infty)}(t)}{\sum_{i=1}^n (X_i-Y_j^*)\indic_{[Y_j^*,\infty)}(X_i)+1/\beta}\right)^{-n_j}.
\end{multline}
\end{corollary} 
By integrating out the vector of latent variables $\bm{Y}$ in \eqref{mom_ext_gamma} we can obtain an estimate of the posterior moments of $\tilde S(t)$. To this end we will use a Gibbs type algorithm whose steps will be described in Section~\ref{algo}.
\section{Approximated inference via moments}\label{sec:approx}
\subsection{Moment-based density approximation and sampling}\label{sec:approx1}

Recovering a probability distribution from the explicit knowledge of its moments 
is a classical problem in probability and statistics that has received great attention in the literature. 
See, e.g., \cite{provost2005moment}, references and motivating applications therein. 
Our specific interest in the problem is motivated by the goal of determining an approximation of the density function of a distribution supported on $[0,1]$ that equals the posterior distribution of a random survival function evaluated at some instant $t$. This is a convenient case since, as the support is a bounded interval, all the moments exist and uniquely characterize the distribution, see \citet{rao1965linear}. 
Moment-based methods for density functions' approximation can be essentially divided into two classes, namely methods that exploit orthogonal polynomial series  \citep{provost2005moment} and maximum entropy methods \citep{csiszar1975divergence,mead1984maximum}. Both these procedures rely on systems of equations that relate the moments of the distribution with the coefficients involved in the approximation. For our purposes the use of orthogonal polynomial series turns out to be more convenient 
since it ensures faster computations as it involves uniquely linear equations. This property is particularly important in our setting since we will need to implement the same approximation procedure for a large number of times in order to approximate the posterior distribution of a random survival function. Moreover, as discussed in \citet{epifani2009moment}, maximum entropy techniques can lead to numerical instability. 

Specifically, we work with Jacobi polynomials, a broad class which includes, among others, Legendre and Chebyshev polynomials. They are well suited for the expansion of densities with compact support contrary to other polynomials like Laguerre and Hermite which can be preferred for densities with infinite of semi-infinite support \citep[see][]{provost2005moment}. 
While the classical Jacobi polynomials are defined on $[-1,1]$, we consider a suitable transformation of such polynomials so that their support coincides with $[0,1]$ and therefore matches the support of the density we aim at approximating. That is, we consider a sequence of polynomials $(G_i)_{i\geq 0}$ such that, for every $i\in \Nbb$, $G_i$ is a polynomial of order $i$ defined by $G_i(s)=\sum_{r=0}^i G_{i,r}s^r$, with $s\in [0,1]$.  The coefficients $G_{i,r}$ can be defined by a recurrence relation %, or as a solution of a differential equation %via the hypergeometric function 
 \citep[see for example][]{szego1967orthogonal}. Such polynomials are orthogonal with respect to the $L^2$-product 
\begin{equation*}\label{eq:scalar_product}
\langle F,G \rangle=\int_0^1 F(s) G(s) w_{a,b}(s)\ddr s,
\end{equation*}
where 
\[w_{a,b}(s)=s^{a-1}(1-s)^{b-1}\]
is named \emph{weight function} of the basis. Moreover, without loss of generality, we assume that the $G_i$'s are normalized and, therefore, $\langle G_i,G_j \rangle=\delta_{ij}$ for every $i,j\in \Nbb$, where $\delta_{ij}$ is the Kronecker symbol. Any univariate density $f$ supported on $[0,1]$ can be uniquely decomposed on such a basis and therefore there is a unique sequence of real numbers $(\lambda_i)_{i \geq 0}$ such that
\begin{equation}\label{eq:decomposition}
f(s)=w_{a,b}(s)\sum_{i=0}^\infty \lambda_i G_i(s).
\end{equation}
Let us now consider a random variable $S$ whose density $f$ has support on $[0,1]$. 
We denote its raw moments by $\mu_r=\E\big[S^r\big]$, with $r\in\Nbb$. 
%$\tilde{S}(t)$ by $\mu_r=\E\big[\big(\tilde{S}(t)\big)^r\big]$. 
From the evaluation of  %$\langle f,G_i\rangle$ 
$\int_0^1 f(s)\, G_i(s)\,\ddr s$ it follows that each $\lambda_i$ coincides with a linear combination of the first $i$ moments, specifically $\lambda_i=\sum_{r=0}^i G_{i,r}\mu_r$. % where the coefficients $G_{i,r}$ are polynomial $G_i$ coefficients. 
Then, the polynomial approximation method consists in truncating the sum in \eqref{eq:decomposition} at a given level $i=N$. This procedure leads to a methodology that makes use only of the first $N$ moments and provides the approximation
\begin{equation}\label{eq:approx}
f_N(s)=w_{a,b}(s)\sum_{i=0}^N \left(\sum_{r=0}^i G_{i,r}\mu_r\right) G_i(s).
\end{equation}
It is important to stress that the polynomial expansion approximation~\eqref{eq:approx} is not necessarily a density as it might fail to be positive or to integrate to 1. In order to overcome this problem, we consider the density $\pi_N$ proportional to the positive part of $f_N$, i.e.  $\pi_N(s)\propto\max(f_N(s),0)$. We resort to importance sampling \citep[see, e.g.,][]{robert2004monte} for sampling from $\pi_N$. 
This is a method for drawing independent weighted samples $(\varpi_\ell, S_\ell)$ from a distribution proportional to a given non-negative function, that exempts us from computing the normalizing constant. More precisely, the method requires to pick a proposal distribution $p$ whose support contains the support of $\pi_N$. A natural choice for $p$ is the Beta distribution %$q\sim\Be(a,b)$ 
proportional to the weight function $w_{a,b}$. % whose support is $(0,1)$.
The weights are then defined by $\varpi_\ell\propto\max(f_N(S_\ell),0)/p(S_\ell)$ such that they add up to 1.

An important issue related to any approximating method refers to the quantification of the approximating error. As for the polynomial approach we undertake, the error can be assessed for large $N$ by applying the asymptotic results in
% the error of approximation can be assessed asymptotically as a function of the number of moments used 
\citet{alexits1961convergence}. In our case, the convergence $f_N(s)\rightarrow f(s)$ for $N\rightarrow \infty$, for all $s\in(0,1)$, implies $\pi_N(s)\rightarrow f(s)$ for $N\rightarrow \infty$. Thus, if $S_N$ denotes a random variable with distribution $\pi_N$, then the following convergence in distribution to the target random variable $S$ holds:
\begin{equation*}%\label{eq:conv-in-distribution}
S_N \stackrel{\mathcal{D}}{\longrightarrow} S \text{ as } N \rightarrow \infty.
\end{equation*}
However, here we are interested in evaluating whether few moments allow for a good approximation of the posterior distribution of $\tilde S(t)$. This question will be addressed by means of an extensive numerical study in the next section. See \citet{epifani2003exponential} and \citet{epifani2009moment} for a similar treatment referring to functionals of neutral-to-the-right priors and Dirichlet processes respectively.

\subsection{Numerical study}\label{sec:approx2}

\begin{center}
\begin{figure}[!ht]
\includegraphics[width=.32\linewidth]{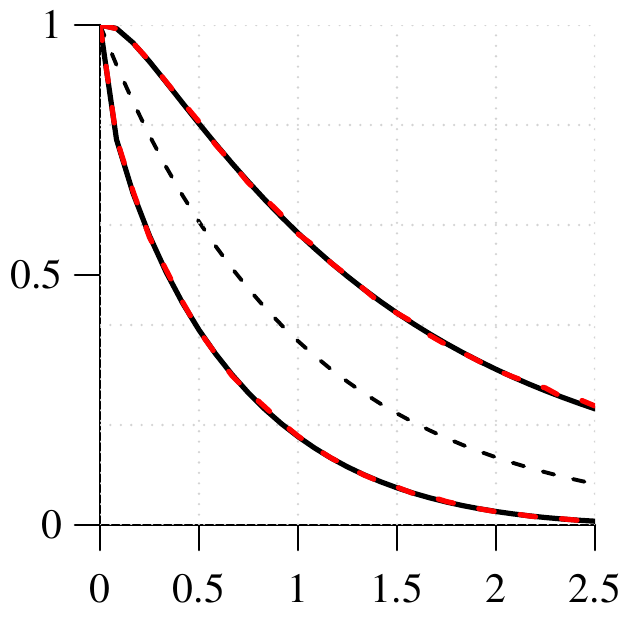}
\includegraphics[width=.32\linewidth]{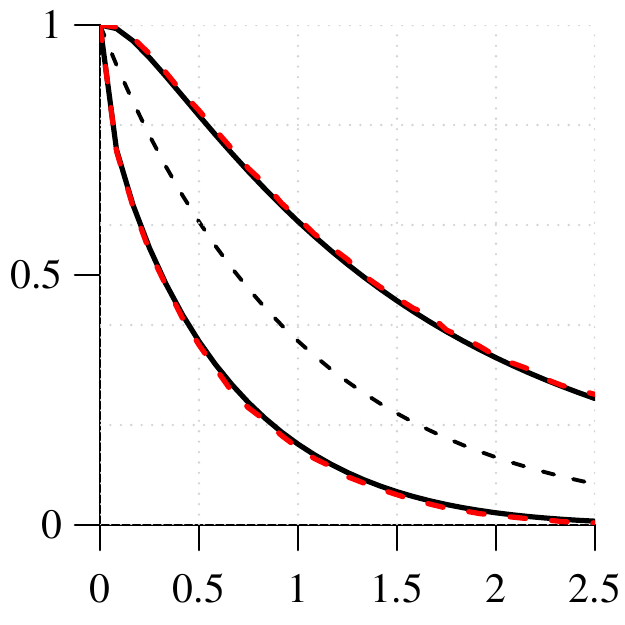}
\includegraphics[width=.32\linewidth]{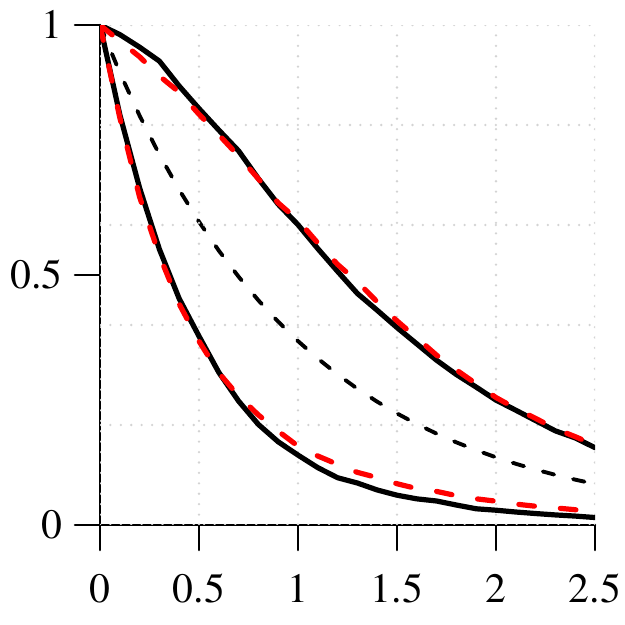}
\caption{Mean of $\tilde S_j(t)$ (dashed black) and $95\%$ highest density intervals for the true distribution (solid black) and the approximated distribution (dashed red) for the Beta ($j=1$), mixture of Beta  ($j=2$) and truncated normal  ($j=3$)  examples (left, middle and right, respectively).}
\label{fig:simulated_examples1}
\end{figure}
\end{center}

\begin{center}
\begin{figure}[!ht]
\includegraphics[width=.98\linewidth]{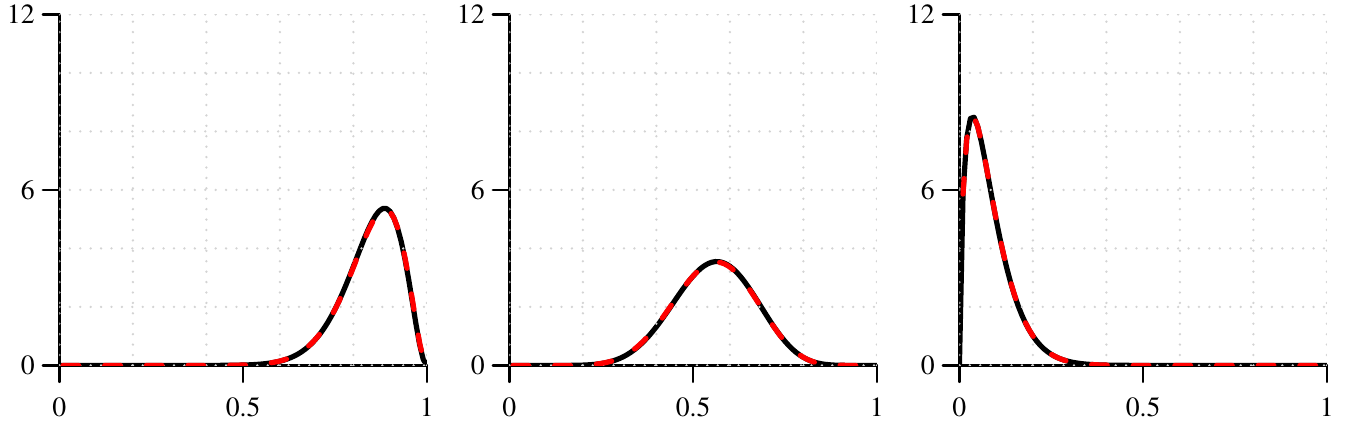}
\includegraphics[width=.98\linewidth]{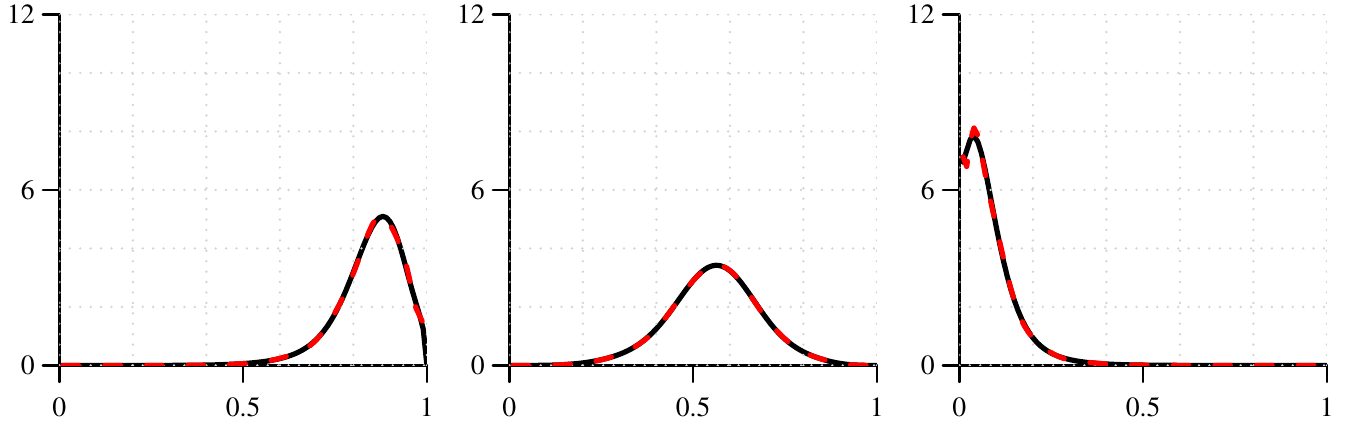}
\includegraphics[width=.98\linewidth]{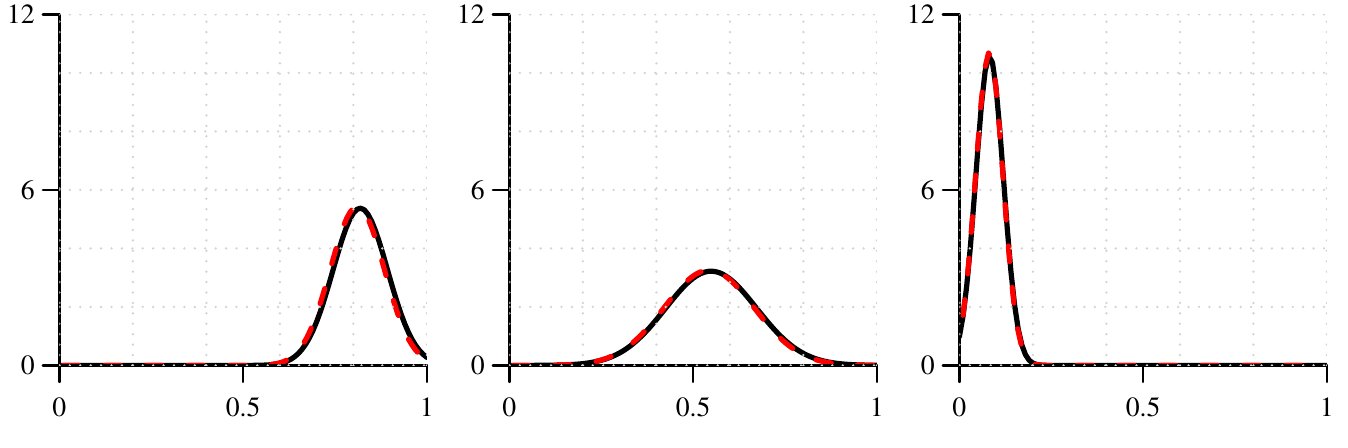}
\caption{True density (solid black) and approximated one (dashed red) at time values $t=0.1$ (left column), $t=0.5$ (middle column) and $t=2.5$ (right column), for the Beta ($j=1$, top row), mixture of Beta ($j=2$, middle row) and truncated normal ($j=3$, bottom row) examples.}
\label{fig:simulated_examples2}
\end{figure}
\end{center}

In this section we assess the quality of the approximation procedure described above by means of a simulation study. The rationale of our analysis consists in considering random survival functions for which moments of any order can be explicitly evaluated at any instant $t$, and then compare the true distribution with the approximation obtained by exploiting the knowledge of the first $N$ moments. This in turn will provide an insight on the impact of $N$ on the approximation error. To this end, we %define 
consider three examples of random survival functions $\tilde S_j$, with $j=1,2,3$. For the illustrative purposes we pursue in this Section, it suffices to specify 
%, by specifying, for every $t> 0$, 
the distribution of the random variable that coincides with $\tilde S_j$ evaluated in $t$, for every $t>0$. Specifically, we consider a Beta, a mixture of Beta, and a normal distribution truncated to $[0,1]$, that is%respectively denoted by $\truncNorm$
\begin{align*}\label{eq:simulated_examples}
\tilde S_1(t)&\sim\Be\left(\frac{S_0(t)}{a_1}, \frac{1-S_0(t)}{a_1}\right), &\nonumber\\[0.3pt]
%&\Unif(S_0(t)+/-aS_0(t)(1-S_0(t))), &\\
\tilde S_2(t)&\sim\frac{1}{2}\Be\left(\frac{S_0(t)}{a_2}, \frac{1-S_0(t)}{a_2}\right)+\frac{1}{2}\Be\left(\frac{S_0(t)}{a_3}, \frac{1-S_0(t)}{a_3}\right), &\\[0.3pt]
\tilde S_3(t)&\sim\truncNorm\left(S_0(t),\frac{S_0(t)(1-S_0(t))}{a_4}\right), & \nonumber
\end{align*}
where %the baseline survival function $S_0$ is defined as an exponential of parameter 1 in all the three cases, i.e. 
$S_0(t)=\edr^{-t}$ and we have set $a_1=20$, $(a_2,a_3)=(10,30)$ and $a_4=2$. 
%In the Beta case, $S_0(t)$ coincides with the mean of the distribution, while it is slightly different for the truncated normal distribution. Note that the $a$ parameter tunes the variance of $\tilde S(t)$. 
We observe that, for every $t>0$, $\E[\tilde S_1(t)]=\E[\tilde S_2(t)]=S_0(t)$ but the same does not hold true for $\tilde S_3(t)$.\\
For each $j=1,2,3$, we computed the first 10 moments of $\tilde S_j(t)$ on a grid $\{t_1,\ldots,t_{50}\}$ of 50 equidistant values of $t$ in the range $[0,2.5]$. The choice of working with 10 moments will be motivated at the end of the section. Then we used the importance sampler described in Section~\ref{sec:approx1} to obtain samples of size 10~000 from the distribution of $\tilde S_j(t_i)$, for each $j=1,2,3$ and $i=1,\ldots,50$. %point of the time discretization. % The results are illustrated on Figures~\ref{fig:simulated_examples1} and~\ref{fig:simulated_examples2}. 
In Figure~\ref{fig:simulated_examples1}, for each $\tilde S_j$, we plot the true mean as well as the $95\%$ highest density intervals for the true distribution and for the approximated distribution obtained by exploiting $10$ moments. %for the Beta, the mixture of Beta and the truncated normal examples (resp. left, middle and right). 
Notice that the focus is not on approximating the mean since moments of any order are the starting point of our procedure. Interestingly, the approximated intervals show a very good fit to the true ones in all the three examples.  As for the Beta case, the fit is exact since the Beta-shaped weight function allows the true density to be recovered with the first two moments. As for the mixture of Beta, exact and approximated intervals can hardly be distinguished. Finally, the fit is pretty good also for the intervals in the truncated normal example.
% in the truncated normal example, the approximated intervals are a little translated compared to the true ones, however they capture most of the true density.
Similarly, in Figure~\ref{fig:simulated_examples2} we compare the true and the approximated densities of each $\tilde S_j(t)$ for fixed $t$ in $\{0.1,0.5,2.5\}$. Again, all the three examples show a very good pointwise fit.
% The pointwise estimate of the density are somewhat rougher than the true one, but still very similar.
$\,$\\[-69pt]

\begin{center}
\begin{figure}[!ht]
\begin{center}
\includegraphics[width=.5\linewidth]{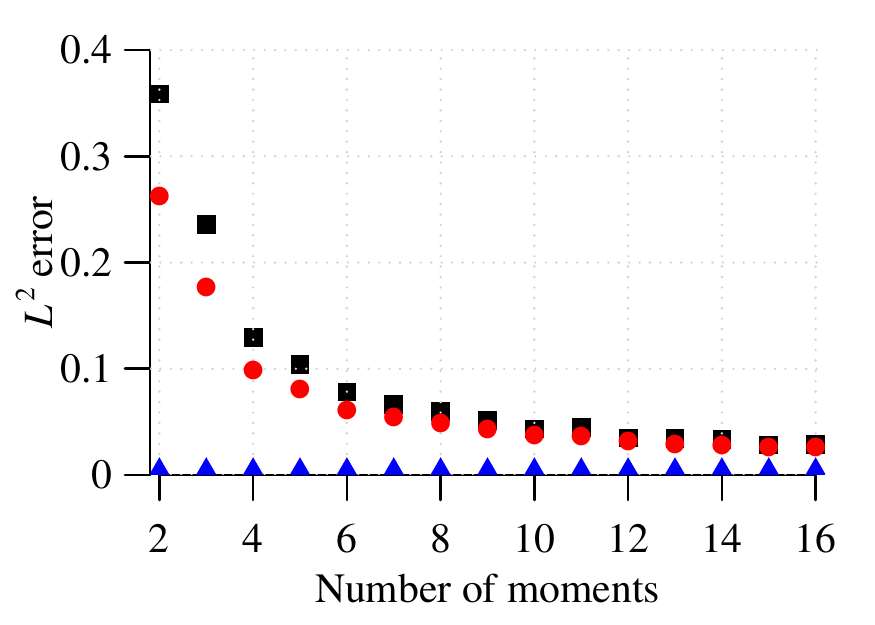}
\caption{Average across $t$ of the $L^2$ error between the true and the approximated densities of $\tilde S_j(t)$, in the Beta example (blue triangles), the mixture of Beta (red dots) and the truncated normal example (black squares). The approximation is exact in the Beta example.}
\label{fig:L2_error}
\end{center}
\end{figure}
\end{center}

We conclude this section by assessing how the choice of $N$ affects the approximation error. To this end, for each instant $t$ on the grid, we numerically compare the true and approximated distributions of $\tilde S_j(t)$, by computing the integrated squared error ($L^2$ error) between the two. Thus, we consider as measure of the overall error of approximation the average of these values. The results are illustrated in Figure~\ref{fig:L2_error}. As expected, the approximation is exact in the Beta example. In the two other cases, we observe that the higher is the number of exploited moments, the lower is the average approximation error. Nonetheless, it is apparent that the incremental gain of using more moments is more substantial when $N$ is small whereas it is less impactful as $N$ increases: for example in the mixture of Beta case, the $L^2$ error is 2.11, 0.97, 0.38 and 0.33 with $N$ equal to 2, 4, 10 and 20 respectively. Moreover, when using a large number of moments, e.g. $N>20$, some numerical instability 
can occur.
These observations suggest that working with $N=10$ moments in~\eqref{eq:approx} strikes a good balance between accuracy of approximation and numerical stability. 
\section{Bayesian inference}\label{sec:real}
In this section we combine the characterization of the posterior moments of $\tilde S(t)$ provided in Proposition~\ref{mom}  together with the approximation procedure described in Section~\ref{sec:approx1}. The model specification \eqref{eq:exch_model} is completed by assuming an extended gamma prior for $\tilde h(t)$, with exponential base measure $P_0(\ddr y)=\lambda\exp(-\lambda y)\ddr y$, and considering the hyperparameters $c$ and $\beta$ random. This leads to the expression \eqref{mom3} for the posterior characterization of the moments. Finally we choose for both $c$ and $\beta$ independent gamma prior distributions with shape parameter 1 and rate parameter $1/3$ (so to ensure large prior variance)  and set $\lambda=1$. Given a sample of survival times $\bm{X}=\{X_1,\ldots,X_n\}$, we estimate the first $N$ moments of the posterior distribution of $\tilde S(t)$, for $t$ on a grid of $q$ equally-spaced points $\{t_1,\ldots,t_{q}\}$ in an interval $[0,M]$, and then we exploit the estimated moments to approximate the posterior distribution of $\tilde S(t_i)$ for $i=1,\ldots,q$. This allows us to devise an algorithm for carrying out full Bayesian inference on survival data. In the illustrations we will focus on the estimation of the median survival time and, at any given $t$ in the grid, of the posterior mean, posterior median, posterior mode and 
credibility intervals for $\tilde S(t)$. The same approach can be, in principle, used to estimate other functionals of interest.

\subsection{Algorithm}\label{algo}
The two main steps we take for drawing samples from the posterior distribution of $\tilde S(t)$, for any $t\in\{t_1,\ldots,t_q\}$, are summarized in Algorithm~\ref{algo:sampler}. First we perform a Gibbs sampler for marginalizing the latent variables $\bm{Y}$ and the hyperparameters $(c,\beta)$ out of \eqref{mom3} and therefore, for every $i=1,\ldots,q$, we obtain an estimate for the posterior moments $\E[\tilde S^{r}(t_i)\vert \bm{X}]$, with $r=1,\ldots,N$. We run the algorithm for $l_{\max}=100~000$ iterations, with a burn-in period of $l_{\min}=10~000$. Visual investigation of the traceplots of the parameters, in the illustrations of Sections~\ref{sec:simulated}~and~\ref{sec:leukemia}, did not reveal any convergence issue. The second part consists in sampling from the posterior distribution of $\tilde S(t_i)$, for every $i=1,\ldots,q$, by means of the importance sampler described in Section~\ref{sec:approx1}. Specifically we sample $\ell_{\max} = 10~000$ values for each $t_i$ on the grid. 

\begin{center}
\begin{algorithm}[h!]
\caption{Posterior sampling\label{algo:sampler}}\medskip
\textbf{Part 1.} \textbf{Gibbs sampler} 
  \begin{algorithmic}[1]
  \STATE set $l=0$ and admissible values for latent variables and hyperparameters,\\i.e. $\{Y_1=Y_1^{(0)},\ldots,Y_n=Y_n^{(0)}\}$, $c=c^{(0)}$ and $\beta=\beta^{(0)}$\\[7pt]
  \STATE while $l< l_{\max}$,  set $l=l+1$, and
  \begin{itemize}
  \item  update $Y_{j}=Y_j^{(l)}$ by means of \eqref{yfc}, for every $j=1,\ldots,n$
  \item update $c=c^{(l)}$ and $\beta=\beta^{(l)}$ by means of \eqref{cfc} and \eqref{betafc}
  \item if $l >  l_{\min}$, compute  
\begin{equation}\label{momgibbs}
 \mu^{(l)}_{r,t}=\E[\tilde S^r(t)\,\vert\,\bm{X},\bm{Y}^{(l)},c^{(l)},\beta^{(l)}]
 \end{equation}
 by means of \eqref{mom3} for each $r=1,\ldots,N$ and for each $t$ in the grid
  \end{itemize}
  \STATE for each $r=1,\ldots,N$ and each $t$ define %\\[7pt]
  $\hat\mu_{r,t}=\frac{1}{l_{\max}-l_{\min}} \sum_{l=l_{\min}+1}^{l_{\max}} \mu_{r,t}^{(l)}$\\[7pt]
  \end{algorithmic}
\textbf{Part 2.} \textbf{Importance sampler} 
 \begin{algorithmic}[1]
  \STATE  for each $t$, use \eqref{eq:approx} and define the approximate posterior density of $\tilde S(t)$ by %\\[7pt]
  $f_{N,t}(\,\cdot\,)=w_{a,b}(\,\cdot\,)\sum_{i=0}^N \left(\sum_{r=0}^i G_{i,r}\hat\mu_{r,t}\right) G_i(\,\cdot\,)$, where $\hat\mu_{0,t}\equiv 1$\\[7pt]
  \STATE  draw a weighted posterior sample $(\varpi_{\ell,t}, S_{\ell,t})_{\ell = 1,\ldots,\ell_{\max}}$ of $\tilde S(t)$, of size $\ell_{\max}$, from $\pi_{N,t}(\,\cdot\,)\propto \max\big(f_{N,t}(\,\cdot\,),0\big)$ by means of the important sampler described in Section~\ref{sec:approx1}\\[7pt]
 \end{algorithmic}
  \end{algorithm}
\end{center}

The drawn samples allow us to approximately evaluate the posterior distribution of $\tilde S(t_i)$, for every $i=1,\ldots,q$. This, in turn, can be exploited to carry out meaningful Bayesian inference (Algorithm~\ref{algo:inference}). As a remarkable example, we consider the median survival time that we denote by $m$. The identity for the cumulative distribution function of $m$
\begin{equation*}%\label{eq:}
\P\left(m\leq t\vert\bm{X}\right) = \P\big(\tilde S(t) \leq 1/2 \vert\bm{X}\big)
\end{equation*}
allows us to evaluate the CDF of $m$ at each time point $t_i$ as $c_i=\P\big(\tilde S(t_i) \leq 1/2 \vert\bm{X}\big)$. Then, we can estimate the median survival time $m$ by means of the following approximation:
\begin{equation}\label{mst}
\hat m=\E_{\bm{X}}[m]=\int_0^\infty \P[m>t\vert\bm{X}]\:\ddr t\approx \frac{M}{q-1}\sum_{i=1}^q(1-c_i)
\end{equation}
where the subscript $\bm{X}$ in $\E_{\bm{X}}[m]$ indicates that the integral is with respect to the distribution of $\tilde S(\cdot)$ conditional to $\bm{X}$. Equivalently, 
\begin{equation}\label{mst2}
\hat m\approx \sum_{i=1}^q t_i (c_{i+1}-c_i),
\end{equation}
with the proviso that $c_{q+1}\equiv 1$. Moreover, the sequence $(c_i)_{i=1}^q$ can be used to devise credible intervals for the median survival time, cf. Part 1 of Algorithm~\ref{algo:inference}. Note that both in \eqref{mst} and in \eqref{mst2} we approximate the integrals on the left-hand-side by means of simple Riemann sums and the quality of such an approximation clearly depends on the choice of $q$ and on $M$. Nonetheless, our investigations suggest that if $q$ is sufficiently large the estimates we obtain are pretty stable and that the choice of $M$ is not crucial since, for $t$ sufficiently large, $\P\big(\tilde S(t) \leq 1/2\vert \bm{X}\big)\approx 0$. Finally, the posterior samples generated by Algorithm~\ref{algo:sampler} can be used to obtain a $t$-by-$t$ estimation of other functionals that convey meaningful information such as the posterior mode and median (together with the posterior mean), cf. Part 2 of Algorithm~\ref{algo:inference}.

\begin{center}
\begin{algorithm}
\caption{Bayesian inference\label{algo:inference}}\medskip
\textbf{Part 1.} \textbf{Median survival time} %for $\tilde S(t) \,\vert\, \bm{X}$
 \begin{algorithmic}[1]
  \STATE use the weighted sample $(\varpi_{\ell,t_i}, S_{\ell,t_i})_{\ell = 1,\ldots,\ell_{\max}}$ to estimate, for each $i=1,\ldots,q$, $c_i=\P(\tilde S(t_i)\leq 1/2 \vert \bm{X})$\\[7pt]
 \STATE plug the $c_i$'s in \eqref{mst2} to obtain $\hat m$\\[7pt]
 \STATE use the sequence $(c_i)_{i=1}^q$ as a proxy for the posterior distribution of $m$ so to devise credible intervals for $\hat m$.\\[7pt]
   \end{algorithmic}
\textbf{Part 2.} \textbf{$t$-by-$t$ functionals} %for $\tilde S(t) \,\vert\, \bm{X}$
 \begin{algorithmic}[1]
  \STATE use the weighted sample $(\varpi_{\ell,t_i}, S_{\ell,t_i})_{\ell = 1,\ldots,\ell_{\max}}$ to estimate, for each $i=1,\ldots,q$, $a_{i}=\inf_{x\in [0,1]}\{\P(\tilde S(t_i)\leq x \vert \bm{X})\geq 1/2\}$ and $b_i=\mbox{mode}\{\tilde S(t_i)\vert \bm{X}\}$\\[7pt]
  \STATE use the sequences $(a_{i})_{i=1}^q$ and $(b_{i})_{i=1}^q$ to approximately evaluate, $t$-by-$t$, posterior median and mode respectively\\[7pt]
  \STATE use 
  the weighted sample $(\varpi_{\ell,t_i}, S_{\ell,t_i})_{\ell = 1,\ldots,\ell_{\max}}$ to devise $t$-by-$t$ credible intervals\\[7pt] 
  \end{algorithmic}
\end{algorithm}
\end{center}

\indent The rest of this section is divided in two parts in which we apply Algorithms~\ref{algo:sampler} and \ref{algo:inference} to simulated and real survival data. In Section~\ref{sec:simulated} we focus on the estimation of the median survival time for simulated samples of varying size. In Section~\ref{sec:leukemia} we analyze a real two-sample dataset and we estimate posterior median and mode, together with credible intervals, of $\tilde S(t)$. In both illustrations our approximations are based on the first $N=10$ moments.

\subsection{Application to simulated survival data}\label{sec:simulated}
\begin{center}
\begin{figure}[ht!]
\includegraphics[width=.47\linewidth]{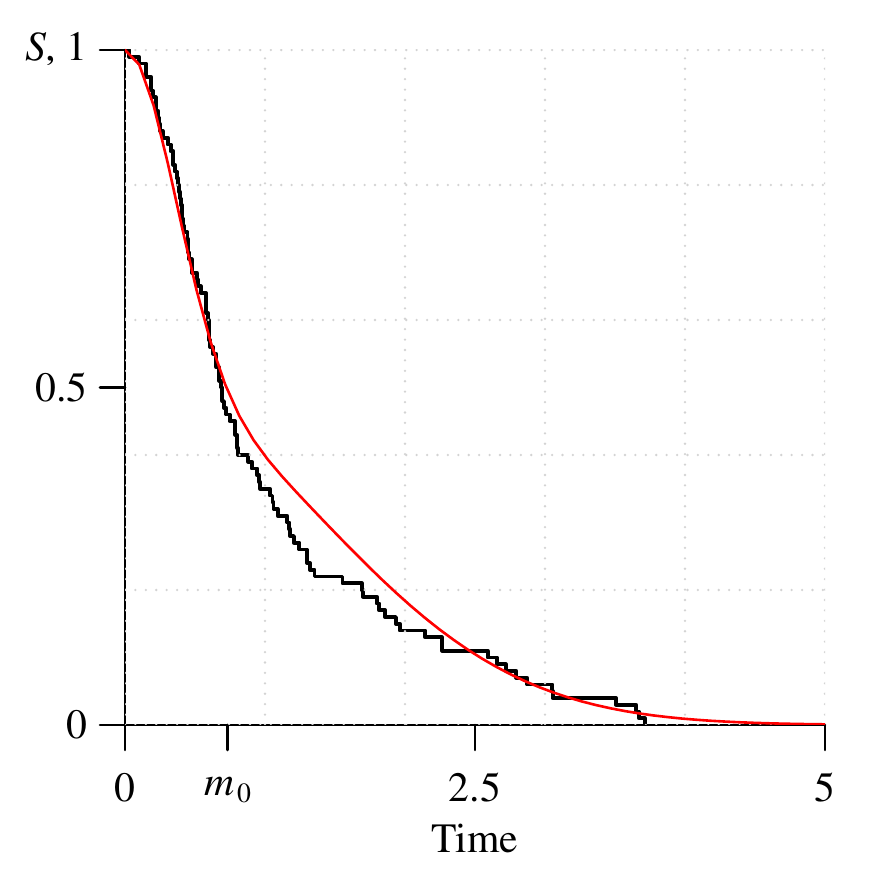}
\includegraphics[width=.51\linewidth]{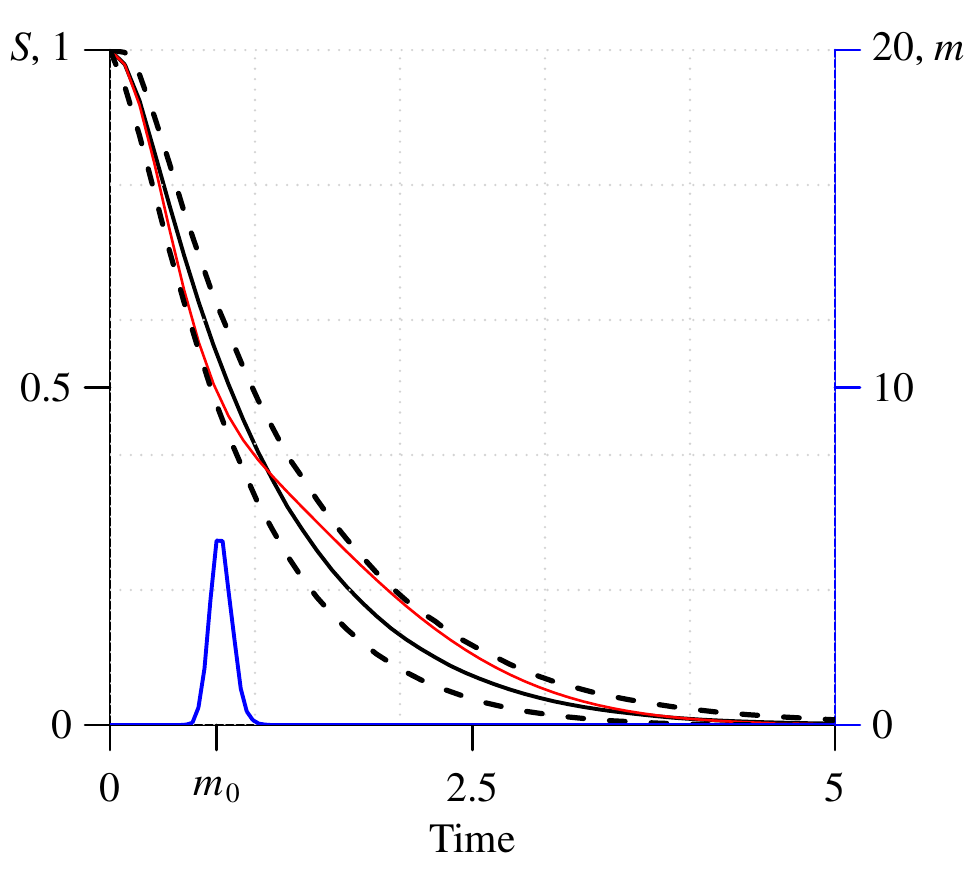}
\caption{(Simulated dataset, $n=100$.) Left: true survival function (red line) and Kaplan--Meier estimate (balk line). Right:  true survival function (red line) and estimated posterior mean (black solid line) with 95\% HPD credible intervals for $\tilde S(t)$ (black dashed lines); the blue plot appearing in the panel on the right is the posterior distribution of the median survival time $m$.}
\label{fig:simulated_data}
\end{figure}
\end{center}

We consider four samples of size $n=25,50,100,200$, from a mixture $f$ of Weibull distributions, defined by
\begin{equation*}
f = \frac{1}{2}\mbox{Wbl}(2,2)+\frac{1}{2}\mbox{Wbl}(2,1/2).
\end{equation*}
After observing that the largest observation in the samples is 4.21, we set $M=5$ and $q=100$ for the analysis of each sample. By applying Algorithms~\ref{algo:sampler} and \ref{algo:inference} we approximately evaluate, $t$-by-$t$, the posterior distribution of $\tilde S(t)$ together with the posterior distribution of the median survival time $m$. In Figure~\ref{fig:simulated_data} we focus on the sample corresponding to $n=100$. On the left panel, true survival function and Kaplan--Meier estimate are plotted. By investigating the right panel we can appreciate that the 
estimated HPD credible regions 
for $\tilde S(t)$ contain the true survival function. Moreover, the posterior distribution of $m$ is nicely concentrated around the true value $m_0=0.724$.
We have investigated the performance of our methodology as the sample size $n$ grows. Table~\ref{table_mst} summarizes the values we obtained for $\hat m$ and the corresponding credible intervals. For all the sample sizes considered, credible intervals for $\hat m$ contain the true value. Moreover, as expected, as $n$ grows, they shrink around $m_0$: for example the length of the interval reduces from 0.526 to 0.227 when the size $n$ changes from 25 to 200. Finally, for all these samples, the estimated median survival time $\hat m$ is closer to $m_0$ than the empirical estimator $\hat m_e$. 
\begin{table}[h!]
\caption{(Simulated datasets.) Comparison of the estimated median survival time ($\hat m$) obtained by means of our Bayesian nonparametric procedure (BNP) and the empirical median survival time $\hat m_e$, for different sample sizes. For BNP estimation we show $\hat m$, the absolute error $|\hat m-m_0|$ and the 95\%-credible interval (CI); last two columns show the empirical estimate $\hat m_e$ and the corresponding absolute error $|\hat m_e - m_0|$. The true median survival time $m_0$ is 0.724.}
\begin{center}
\begin{tabular}{ l l l l l l l l }
  \hline
&  &\multicolumn{3}{c}{BNP}  & & \multicolumn{2}{c}{Empirical}\\
\hline
 sample size & & $\hat m$ & error & CI& & %C.I. length &  
  $\hat m_e$ & error\\ 
  \hline
   25 && 0.803 & 0.079 & (0.598, 1.124)% & 0.527 \\ 
& & 0.578 & 0.146 \\ 
  50 && 0.734 & 0.010 & (0.577, 0.967) %& 0.389 \\ 
& &0.605   & 0.119 \\ 
  100& & 0.750 & 0.026 & (0.622, 0.912)% & 0.290 \\ 
& & 0.690  & 0.034 \\ 
  200 && 0.746 & 0.022 & (0.669, 0.896) %& 0.227
& & 0.701   & 0.023 \\ 
   \hline
\end{tabular}\label{table_mst}
   \end{center}
   \end{table}

\subsection{Application to real survival data}\label{sec:leukemia}
\begin{figure}[ht!]
\begin{center}
\includegraphics[width=.323\linewidth]{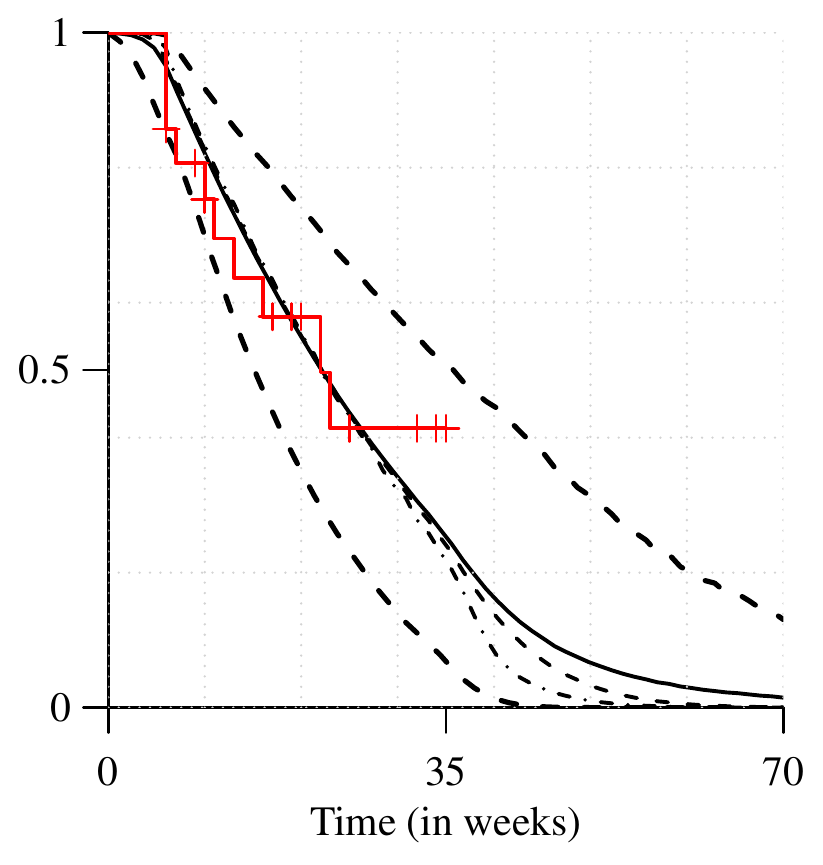}
\includegraphics[width=.323\linewidth]{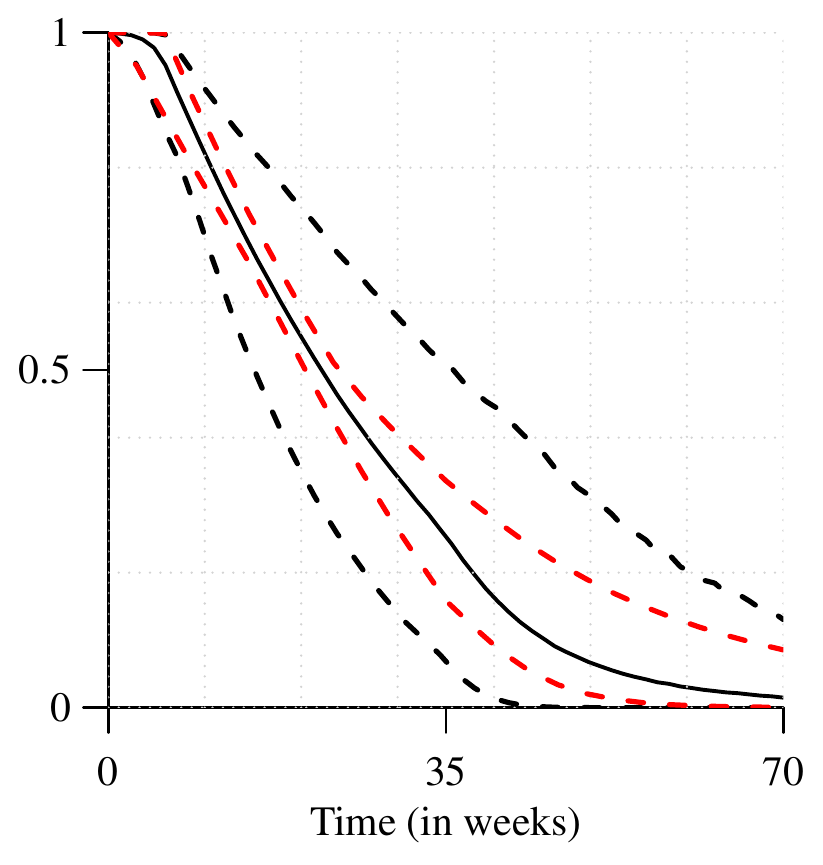}
\includegraphics[width=.323\linewidth]{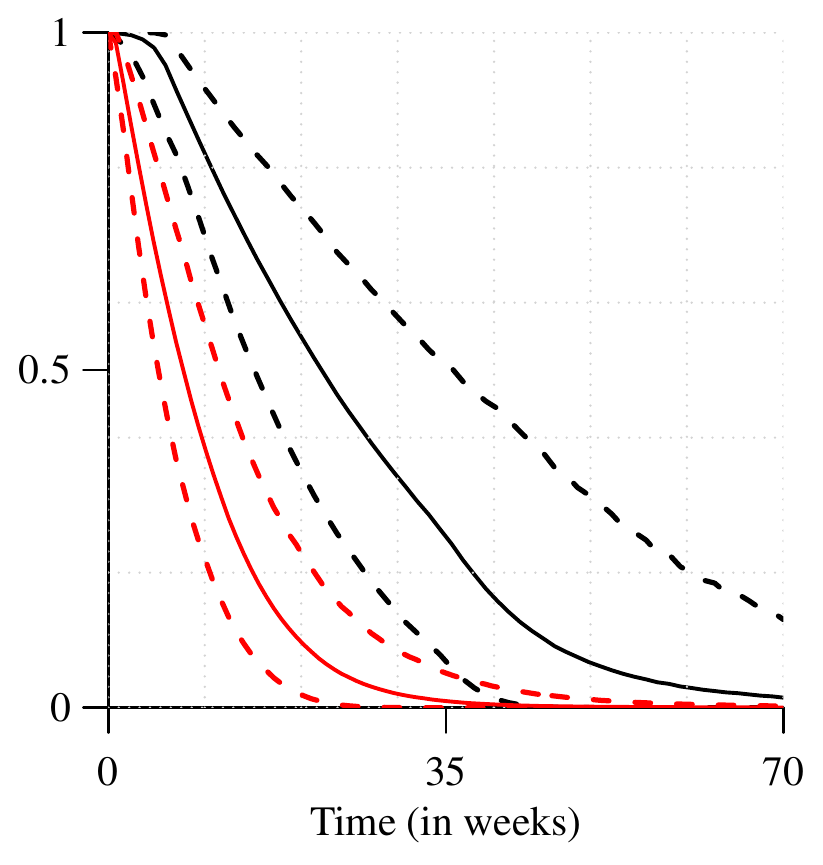}
\caption{Left: comparison of posterior mean (solid line), median (dashed line) and mode (point dashed line) in dataset $\mathsf{T}$, with 
95\% HPD credible intervals 
(dashed line). The Kaplan--Meier estimate is plotted in red. Middle: comparison of the 
95\% HPD credible interval 
(dashed black line) with the marginal interval (dashed red line).
Right: comparison of samples $\mathsf{T}$ (black) and $\mathsf{P}$ (red), with posterior means (solid) and 
95\% HPD credible intervals 
(dashed).}
\label{fig:simulated_examples}
\end{center}
\end{figure}

We now analyze, with the described methodology, a well known two-sample dataset involving leukemia remission times, in weeks, for two groups of patients, under active drug treatment and placebo respectively. The same dataset was studied, e.g., by \cite{Cox72}.
Observed remission times for patients under treatment ($\mathsf{T}$) are
$$\{6,6,6,6^*,7,9^*,10,10^*
,11,13,16,17^*
,19^*
,20^*
,22,23,25^*
,32^*
,32^*
,34^*
,35^*\},$$
where stars denote right-censored observations. Details on the censoring mechanism and on how to adapt our methodology to right-censored observations are provided in \ref{sec:censored}. On the other side, remission times of patients under placebo ($\mathsf{P}$) are all exact and coincide with 
$$\{1,1,2,2,3,4,4,5,5,8,8,8,11,11,12,12,15,17,22,23\}.$$
For this illustration we set $M=2\max(\bm{X})$, that is $M=70$, and $q=50$. For both samples we estimate and compare posterior mean, median and mode as well as 95\% credible intervals. 
In  the left panel of Figure~\ref{fig:simulated_examples} we have plotted such estimates for sample $\mathsf{T}$. By inspecting the plot, it is apparent that,  for large values of $t$, posterior mean, median and mode show significantly different behaviors, with posterior mean being more optimistic than posterior median and mode. It is worth stressing that such differences, while very meaningful for clinicians, could not be captured by marginal methods for which only the posterior mean would be available. A fair analysis must take into account the fact that, up to $t=23$, i.e. the value corresponding to the largest non-censored observation, the three curves are hardly distinguishable. The different patterns for larger $t$ might  therefore depend on the prior specification of the model. Nonetheless, we believe this example is meaningful as it shows that a more complete posterior analysis is able to capture differences, if any, between posterior mean, median and mode.

When relying on marginal methods, the most natural choice for estimating the uncertainty of posterior estimates consists in considering the quantiles intervals corresponding to the output of the Gibbs sampler, that we refer to as \emph{marginal  intervals}. This leads to consider, for any fixed $t$, the interval whose lower and upper extremes are the quantiles of order $0.025$ and $0.975$, respectively, of the sample of conditional moments
$\{\mu_{1,t}^{(l_{\min}+1)},\ldots,\mu_{1,t}^{(l_{\max})}\}$ 
defined in \eqref{momgibbs}. In the middle panel of Figure~\ref{fig:simulated_examples} we have compared the estimated 95\% HPD intervals for $\tilde S(t)$ and the marginal  intervals corresponding to the output of the Gibbs sampler. In this example, the marginal method clearly underestimates the uncertainty associated to the posterior estimates. 
This can be explained by observing that, since the underlying completely random measure has already been marginalized out, the intervals arising from the Gibbs sampler output, capture only the variability of the posterior mean that can be traced back to the latent variables $\bm{Y}$ and the parameters $(c,\beta)$. As a result, the uncertainty detected by the marginal method leads to credible intervals that can be significantly narrower than the actual posterior credible intervals that we approximate through the moment-based approach. This suggests that the use of intervals produced by marginal methods as proxies for posterior credible intervals should be, in general, avoided.

We conclude our analysis by observing that the availability of credible intervals for survival functions can be of great help in comparing treatments. In the right panel of Figure~\ref{fig:simulated_examples} posterior means as well as corresponding 95\% HPD intervals are plotted for both samples $\mathsf{T}$ and $\mathsf{P}$. By inspecting the plot, for example, the effectiveness of the treatment seems clearly significant as, essentially, there is no overlap between credible intervals of the two groups.

\section*{Acknowledgment}
J. Arbel and A. Lijoi are supported by the European Research Council (ERC) through StG ``N-BNP'' 306406.

\appendix

\section{Moments under exponential $P_0$}\label{app_mom} 
We provide an explicit expression for \eqref{mom_ext_gamma} when $P_0(\ddr y)=\lambda \exp(-\lambda y )\ddr y$ and the hyperparameters $c$ and $\beta$ are considered random.
\begin{multline}\label{mom3}
\E[\tilde S^r(t)\,|\,\bm{X},\bm{Y},c,\beta]=\\
\exp\left\{-c\edr^{-f_{0,r}(0)}\left[\mbox{Ei}(f_{0,r}(t))-\mbox{Ei}(f_{0,r}(X_1\wedge t))\right]\right\}\left(\frac{f_{0,r}(X_1\wedge t)}{f_{0,r}(t)}\right)^{-c\edr^{-\lambda (X_1\wedge t)}}\\
\times\prod_{i=1}^n\exp\left\{-c\edr^{-f_{i,r}(0)}\left[\mbox{Ei}(f_{i,r}(X_i\wedge t))-\mbox{Ei}(f_{i,r}(X_{i+1}\wedge t))\right]\right.\\
\left. -c\edr^{-f_{i,0}(0)}\left[\mbox{Ei}(f_{i,0}(X_{i+1}\wedge t))-\mbox{Ei}(f_{i,0}(X_{i}\wedge t))\right]\right\}\\[5pt]
\times\left(\frac{i}{i+r}\frac{f_{i,0}(X_i\wedge t)}{f_{i,r}(X_i\wedge t)}\right)^{-c\edr^{-\lambda (X_i\wedge t)}}\left(\frac{i+r}{i}\frac{f_{i,r}(X_{i+1}\wedge t)}{f_{i,0}(X_{i+1}\wedge t)}\right)^{-c\edr^{-\lambda (X_{i+1}\wedge t)}}\\
\times\prod_{j=1}^k \left(1+r\,\frac{(t-Y_j^*)\indic_{[Y_j^*,\infty)}(t)}{\sum_{i=1}^n (X_i-Y_j^*)\indic_{[Y_j^*,\infty)}(X_i)+1/\beta}\right)^{-n_j},
\end{multline} 
where $\mbox{Ei}(\cdot)$ is the exponential integral function defined for non-zero real values $z$ by 
$$\mbox{Ei}(z)=-\int_{-z}^\infty\frac{\edr^{-t}}{t}\ddr t$$
 and the function $f_{i,r}$, for $i,r\geq 0$ such that $i+r>0$, is defined by
$$f_{i,r}(x)=\lambda\left(\frac{\xi_{i}+1/\beta+rt}{i+r}-x\right).$$

\section{Full conditional distributions}\label{app_fc}
In this section we provide expressions for the full conditional distributions needed in the algorithm described in Section~\ref{algo} for extended gamma processes with base measure $P_0(\ddr y)=\lambda \exp(-\lambda y) \ddr y$. These distributions are easily derived, up to a constant, from the joint distribution of the vector $(\bm{X},\bm{Y},c,\beta)$, that can be obtained from \eqref{eq:joint.rhr.1}.
Therefore we start by providing the full conditional distribution for the latent variable $Y_i$, with $i=1,\ldots,n$, where $\bm{Y}^{(-i)}$ denotes the vector of distinct values $(\tilde{Y}^*_1,\ldots,\tilde{Y}^*_{k^*})$ in $(Y_1,\ldots,Y_{i-1},Y_{i+1},\ldots,Y_n)$ and $(n_1^{(-i)},\ldots,n_{k^{*}}^{(-i)})$ represent the corresponding frequencies.
\begin{equation}\label{yfc}
\P[Y_i = \ddr y \, | \, \bm{X}, \bm{Y}^{(-i)},c,\beta]=p_0 G_0 (\ddr y)+ \sum_{j=1}^{k^*} p_j \delta_{\tilde{Y}_j^{*}}(\ddr y),
\end{equation}
where
\begin{align*}
p_0&\propto c\,\lambda\, \sum_{j=i}^n \frac{1}{j}\edr^{-\lambda \frac{\xi_j+1/\beta}{j}}\left[\mbox{Ei}\left(f_{j,0}(X_{j+1})\right)-\mbox{Ei}\left(f_{j,0}(X_j)\right)\right],\\
p_j&\propto\indic_{\{Y_j^*\leq X_{i}\}}\frac{n_j^{(-i)}}{\sum_{l=1}^n (X_l-\tilde Y_j^*)\indic_{[0,X_l)}(\tilde Y_j^*)+1/\beta}
\end{align*} 
and
\begin{equation*}
G_0(\ddr y)\propto \indic_{[0,X_i)}(y)\edr^{-\lambda y}\frac{1}{\sum_{j=1}^n (X_j-y)\indic_{[0,X_j)}(y)+1/\beta}\ddr y.
\end{equation*}
Finally, the full conditional distributions for the parameters $c$ and $\beta$ are given respectively by
\begin{multline}\label{cfc}
\mathcal{L}(c\,|\,\bm{X},\bm{Y},\beta)\propto \mathcal{L}_{0}(c) c^k \beta^{-c}\prod_{i=1}^n \exp\left\{-c \edr^{-f_{i,0}(0)}\left[\mbox{Ei}(f_{i,0}(X_i))-\mbox{Ei}(f_{i,0}(X_{i+1}))\right] \right\}\\
\times \frac{\left(\xi_i+1/\beta-iX_{i+1}\right)^{-c\edr^{-\lambda X_{i+1}}}}{\left(\xi_i+1/\beta-iX_{i}\right)^{-c\edr^{-\lambda X_i}}}
\end{multline}
and
\begin{multline}\label{betafc}
\mathcal{L}(\beta\,|\,\bm{X},\bm{Y},c)\propto \mathcal{L}_0(\beta)\beta^{-c}\prod_{i=1}^n \exp\left\{-c \edr^{-f_{i,0}(0)}\left[\mbox{Ei}(f_{i,0}(X_i))-\mbox{Ei}(f_{i,0}(X_{i+1}))\right] \right\}\\
\times \frac{\left(\xi_i+1/\beta-iX_{i+1}\right)^{-c\edr^{-\lambda X_{i+1}}}}{\left(\xi_i+1/\beta-iX_{i}\right)^{-c\edr^{-\lambda X_i}}}\prod_{j=1}^k \left(\sum_{i=1}^n(X_i-Y_j^*)\mathds{1}_{[Y_j^*,\infty)}(X_i)+1/\beta\right)^{-n_j},
\end{multline}
where $\mathcal{L}_0(c)$ and $\mathcal{L}_0(\beta)$ are the prior distributions of $c$ and $\beta$ respectively.
\section{Censored observations}\label{sec:censored}

The methodology we have presented in Section~\ref{sec:real} needs to be adapted to the presence of right-censored observations in order to be applied to the dataset in Section~\ref{sec:leukemia}. Here we introduce some notation and illustrate how the posterior characterization of Proposition~\ref{Jam} changes when data are censored. To this end, let $C_i$ be the right-censoring time corresponding to $X_i$, and define $\Delta_i = \indic_{(0,C_i]}(X_i)$, so that $\Delta_i$ is either 0 or 1 according as to whether $X_i$ is censored or exact. The actual $i$th observation is $T_i=\min(X_i,C_i)$ and, therefore, data consist of pairs $\bm{D} = \{(T_i,\Delta_i)\}_{i=1\ldots n}$.
In this setting, the likelihood in~\eqref{eq:likelihood} can be rewritten as
\begin{equation*}%\label{eq:likelihood}
 \Lc(\mt;\bm{D})=e^{-\int_\Y K^*_{\bm{D}}(y)\mt(\ddr y)}\prod_{i:\,\Delta_i=1}\int_\Y k(T_i;y)\mt(\ddr y),
\end{equation*}
where $$K^*_{\bm{D}}(y)=\sum_{i=1}^n\int_0^{T_i}k(s;y)\ddr s.$$
By observing that the censored times are involved only through $K^*_{\bm{D}}$, we have that the results derived in Proposition~\ref{Jam} under the assumption of exact data easily carry over to the case with right-censored data. The only changes refer to $K_{\bm{X}}$, that is replaced by $K^*_{\bm{D}}$, and the jump components which occur only at the distinct values of the latent variables that correspond to exact observations. For instance in Proposition~\ref{Jam}, the L\'evy intensity of the part of the CRM without fixed points of discontinuity is modified by 
\begin{equation*}
 \nu^*(\ddr s, \ddr y)= e^{-s K_{\bm{D}}^*(y)}\rho_y(s)\ddr s\,c P_0(\ddr y),
\end{equation*}
while the distribution of the jump $J_j$ has density function  $f(\,\cdot\,|\, n_j^*,K_{\bm{D}}^*(Y_j^*),Y_j^*)$ with $f$ defined in \eqref{eq:jump_dens} and $n_j^* = \#\big\{i\,:\,Y_i=Y_j^* \text{ and } \Delta_i = 1\big\}$. Adapting the results of Proposition~\ref{mom} and Corollary~\ref{cor:moments}, as well as the full conditional distributions in \ref{app_fc}, is then straightforward.

\bibliographystyle{apalike}
\bibliography{hazbib}
\end{document}